\documentclass[10pt]{article}
\usepackage{amsmath,amssymb,amsthm,amscd}
\numberwithin{equation}{section}

\def\p{\partial}

\newtheorem{prop}{Proposition}[section]
\newtheorem{theo}[prop]{Theorem}
\newtheorem{lem}[prop]{Lemma}
\newtheorem{cor}[prop]{Corollary}
\newtheorem{rem}[prop]{Remark}

\newtheorem{defi}[prop]{Definition}

\newtheorem{q}[prop]{Question}
\def\begeq{\begin{equation}}
\def\endeq{\end{equation}}

\def\and{\quad{\rm and}\quad}

\let\lra=\longrightarrow

\def\mapright\#1{\,\smash{\mathop{\lra}\limits^{\#1}}\,}

\begin{document}

\makeatletter      
\renewcommand{\ps@plain}{%
     \renewcommand{\@oddhead}{\textrm{Space of K\"{a}hler metrics}\hfil\textrm{\thepage}}%
     \renewcommand{\@evenhead}{\@oddhead}%
     \renewcommand{\@oddfoot}{}
     \renewcommand{\@evenfoot}{\@oddfoot}}
     \renewcommand{\thefootnote}{\fnsymbol{footnote}}
\makeatother     


\title{The Space of K\"{a}hler metrics (II)}
\author{ E. Calabi and X. X. Chen\footnote{The second author is supported
partially by NSF postdoctoral fellowship.}
  }
\date{}
\pagestyle{plain}
\bibliographystyle{plain}
\maketitle
\section{Introduction and Main results}
This paper, the second of a series, deals with the function space
of all smooth K\"ahler metrics in any given closed complex
manifold $M$ in a fixed cohomology class. This function space is
equipped with a pre-Hilbert manifold structure introduced by T.
Mabuchi \cite{Ma87}, where he also showed formally it has
non-positive curvature. The previous result of the second author
\cite{chen991} showed that the space is a path length space and
it is geodesically convex in the sense that any two points are
joined by a unique path, which is always length minimizing and of
class $C^{1,1}.\;$ This already confirms one of Donaldson's
conjecture completely and verifies another one partially (cf.
\cite{Dona96}). In the present paper, we show  first of all, that
the space is, as expected, a path length space of non-positive
curvature in the sense of A. D. Alexanderov.  A second result is
related to the theory of extremal K\"ahler metrics, namely that
the gradient flow of the K energy\footnote{The K energy is
defined by T. Mabuchi in 1987 \cite{Ma87} while the flow was
introduced by the first author in 1982 \cite{calabi82}. And it is
commonly known as the ``Calabi flow" in the literature.} is
strictly length decreasing on all paths except  those induced by
a path of holomorphic automorphisms of $M.\;$ This result, in
particular, implies that extremal K\"ahler metric is unique up to
holomorphic transformations, provided that Donaldson's conjecture
on the regularity of geodesic is true.

\subsection{ Riemannian metrics and Non Positive Curved
space.} Let $(V,\omega_0)$ be a polarized K\"ahler  manifold
without boundary.  Consider the space of K\"ahler distortion
potentials
\[
 {\cal{H}} = \{ \varphi \in C^{\infty}(V) : \omega_{\varphi} =  \omega_{0} + \partial \overline{\partial} \varphi > 0 \; {\rm on}\; V\}.
\]

Clearly, the tangent space  $T \cal {H} $  is $C^{\infty} (V).\;$
Each K\"{a}hler potential  $\phi \in \cal {H}$ defines a measure
$d\,\mu_{\phi} = {1\over {n!}} \omega_{\phi}^n.\; $
 A Weil-Peterson type metric was defined on this infinite
dimensional manifold $\cal {H},$  using the $L^2 $ norm provided
 by these measures (cf. Section 2.1  for historical remark
of this metric). In 1997, following a program of Donaldson
\cite{Dona96},  the second author proved that this space is convex
by $C^{1,1}$ geodesic; and used this fact
 to prove that $\cal H$ is indeed a path length space. If one could prove that the resulting geodesic
 were $C^4$ instead of $C^{1,1}$, then the formal calculations in
\cite{Ma87} would yield that the curvature is non-positive. Since
we only have $C^{1,1}$  geodesics,  an important question is
whether this is a non positive curved space in the sense of A.D.
Aleksandrov. We give an affirmative  answer to this question here:

\begin{theo} The space of K\"ahler  potentials
 in a fixed
K\"ahler  class is Non-Positive Curved space: Suppose A, B, C are
three points in the space of K\"ahler  potentials and
$P_{\lambda}$
 is a geodesic interpolation point of $B$ and $C$ for $ 0 \leq \lambda \leq 1:$
 the geodesic distance from $P_\lambda$ to $B,$ and the distance
 from $P_\lambda$ to $C$ are $\lambda d(B,C)$ and $(1-\lambda)
 d(B,C)$ respectively (here we use $d(P,Q)$ to denote the distance between $P$ and $Q$ in ${\cal H}.$).
Then the following inequality holds:

\[
d(A, P_{\lambda})^2 \leq (1-\lambda) d(A, B)^2 + \lambda d (A, C)^2 - \lambda \cdot (1-\lambda) d(B,C)^2.
\]
\end{theo}

In \cite{chen991}, the second author proved that the geodesic
minimizes all possible length. However,  it is not clear whether a
sequence of curves which minimize the length between any two
points in $\cal H$ converges to a geodesic or not. Now we can give
an affirmative answer to this question:

\begin{theo} For any two metrics $\varphi_0, \varphi_1$ in $\cal H,$
let $\{C_i\}$ be any sequence of curves in $\cal H$ which connect
between $\varphi_0$ and $\varphi_1.\;$ Suppose the length of this
sequence of curves aprroaches  the infimum of length over all
possible curves between $\varphi_0$ and $\varphi_1,$ then $C_i$
converges to the unique $C^{1,1}$ geodesic which connects
$\varphi_0$ and $\varphi_1 \;$ in the sense of distance.
\end{theo}
\subsection{The gradient flow}
 In \cite{calabi82}, the first author
 introduced the notion of extremal K\"ahler metrics and
proposed to use  a fourth order heat equation to
attack the existence problem of extremal K\"ahler  metrics:
\begin{equation}
  {{\partial \, \varphi}\over{\partial\,s}} = R (\varphi) -
  \underline{R},
\label{eq:calabiflow}
\end{equation}
where $R(\varphi)$ is the scalar curvature of the K\"ahler metric
$\omega_{\varphi}$ and $\underline{R}$ is the average scalar
curvature --- a constant depending only the K\"ahler class.
Somewhat surprisingly,   we observed that this flow actually
decreases the length of any smooth curve in $\cal H:$

\begin{theo}  Given
 any two K\"ahler potentials  $\varphi_1$ and
$\varphi_2$ in $\;\cal H$ and a smooth curve $C(t)$ in $\cal H$
 connecting them,
the length of this curve  strictly decreases under gradient flow
(\ref{eq:calabiflow}) unless this curve in $\cal H$ represents a
path of holomorphic transformations.
 More specifically, if $\varphi(t), 0 \leq t \leq 1 $ is
a curve in $\cal H,$ and $L$ is the length of this curve; and
suppose $\varphi(s,t) $ is  the family of curves under the
gradient flow  (\ref{eq:calabiflow}),  then
\[
  {{d \, L}\over {d\, s}} = -  \int_0^1 \left(\int_V |D
   {{\partial \varphi} \over {\partial t}}|_{\varphi(s,t)}^2 \; d \,g(s,t) \cdot \sqrt{ \int_V |{{\partial \varphi} \over {\partial t}}|^2\; d \,g(s,t)}^{-{1\over 2}} \right)\; d\,t .\]
 Here $ g(s,t)$ is the K\"ahler metric
corresponding to the K\"ahler potential $\varphi(s,t),\;$ while
$D$ is the Lichernowicz operator (cf. Definition 1.4 below).
\end{theo}

\begin{defi} The Lichernowicz operator $D:$  For any smooth function
$f$ in $V,\; D(f) = \displaystyle \sum_{\alpha,
\beta=1}^n\;f_{,\alpha \beta} dz^{\alpha} \otimes d z^{\beta}\;$
where $f_{,\alpha \beta}$ is the second covariant derivatives of
$f.\;$
\end{defi}

This flow is known to have short time existence. In Riemann
surface, Chrusciel \cite{Chru91} proved that the global existence
and convergence of this flow if there is a constant scalar
curvature metric apriori.  The second author gave a new geometric
proof to the Chrusciel's theorem \cite{chen994}. Following
approach taken in \cite{chen994}, M. Struwe \cite{Struwe2000} gave
a unified treatment of both Ricci flow and Calabi flow in Riemann
surfaces. In higher dimensional  K\"ahler manifold cases, very
little is known about the global existence of the flow.
\begin{theo} The following two statements hold
\begin{enumerate}
\item If the gradient flow  (\ref{eq:calabiflow}) exists for all time for any
smooth initial data, then the distance between any two metrics
decreases under the gradient flow  (\ref{eq:calabiflow}).
 \item If the K energy is weakly convex\footnote{A function $f(t)
 (0\leq t \leq 1)$ is weakly convex if for any $t$, we have $f(t) \leq (1-t)\; f(0) + t\; f(1).\;$
  For this theorem, we additionally assume that $f$ is differentiable at both end points, i.e., $f'(1) \geq f'(0).\;$ } along geodesic,
then the gradient flow  decreases distance. In particular, if the
first Chern class $C_1(V) \leq 0,$ then the gradient flow
(\ref{eq:calabiflow}) decreases distance.
\end{enumerate}
\end{theo}

The following question is very interesting:
\begin{q} If $C_1(V) < 0,$ is the distance function in $\cal H$
strictly decrease under the gradient flow  (\ref{eq:calabiflow})?
\end{q}




{\bf Acknowledgement.} The second author would like to thank
Professor Donaldson for many stimulating discussions. He also
wants to thank
 Professors L. Simon and R. Schoen for their encouragement throughout
this project. He also wants to thank Guofang Wang for reading
through an earlier version of this paper carefully.
\section{${\cal H}$ is a non positive curved space}
In this section, we want to show that $\cal H$ is a non positive
curved space in the sense of Aleksandrov.
\subsection{ A Riemannian metric in the infinite dimensional space.}
 Mabuchi (\cite{Ma87})
 in 1987 defined a Riemannian metric  on the space of K\"ahler metrics,
under which it  becomes (formally) a non-positive curved infinite dimensional
 symmetric space. Apparently unaware of Mabuchi's work,
Semmes \cite{Semmes92}  and Donaldson \cite{Dona96}  re-introduced this same metric again from different
angles. Let us now introduce this metric here.
A tangent vector in $\cal {H}$ is just a function in $V.\;$
For any  vector $\psi \in T_{\varphi} \cal {H}, $ we define the length of this
vector as
\[
\|\psi\|^2_{\varphi} =\int_{V}\psi^2\;d\;\mu_{\varphi}.
\]
For two ``vectors'' $f_1, f_2$ in $T_{\varphi} \cal H,$ we use the standard
notation in Riemannian geometry to denote their inner product:
\[
  \langle f_1,f_2\rangle_{\varphi} = \displaystyle\;\int_V\; f_1\cdot f_2 \;d\;\mu_{\varphi}.
\]
When no confusion is arisen, we just write
\[
  \langle f_1,f_2\rangle = \displaystyle\;\int_V\; f_1\cdot f_2 \;d\;\mu_{\varphi}.
\]
For a path $\varphi(t) \in {\cal {H}} (0\leq t \leq 1),$
the length  is given by
\[
  \int_0^1 \langle {{\p  \varphi} \over {\p t}}(t), {{\p  \varphi} \over {\p t}}(t)\rangle_{\varphi(t)} \; d\,t=  \int_0^1 \; \sqrt{\int_V {{{\p  \varphi} \over {\p t}}}^2 d\,\mu_{\varphi(t)}} \; d\,t
\]
and the geodesic equation is
\begin{equation}
 {{\p^2  \varphi} \over {\p t^2}} - {1\over 2} |\nabla {{\p  \varphi} \over {\p t}}(t)|^2_{ \varphi(t)} = 0,
\label{geodesic}
\end{equation}
where the derivatives and norm in the second term of the left hand
side
are taken with respect to the metric $\omega_{\varphi(t)}.\;$\\

This geodesic equation shows us how to define a connection on the
tangent bundle of ${\cal H}$. The notation is simple if one
thinks of such a connection as a way of differentiating vector
fields along paths. Thus, if $\phi(t)$ is any path in ${\cal H}$
and $\psi(t)$ is a field of tangent vectors along the path (that
is, a function on $V \times [0,1]$), we define the covariant
derivative along the path to be
$$ D_{t}\psi = \frac{\partial\psi}{\partial t} - {1\over 2}  (\nabla \psi, \nabla \phi')_{\phi}.  $$
This connection is torsion-free because in the canonical \lq\lq co-ordinate
chart'', which represents ${\cal H}$ as an open subset of $C^{\infty}(V)$, the
`` Christoffel symbol''
$$\Gamma: C^{\infty}(V) \times C^{\infty}(V) \rightarrow C^{\infty}(V)$$ at
$\phi$ is just
$$\Gamma(\psi_{1}, \psi_{2}) = - {1\over 2}  (\nabla \psi_{1},
\nabla\psi_{2})_{\phi}$$ which is symmetric in
$\psi_{1},\psi_{2}$. It is easy to verify that the connection is
metric-compatible. By a direct calculation,  it was proved
formally in \cite{Ma87} (and later re-proved in \cite{Semmes92}
and
 \cite{Dona96}) that $\cal H$ is a non-positive curved space.

Donaldson \cite{Dona96} in 1996  introduced a  connection between
this formal Riemannian metric in the infinite dimensional space
$\cal H$ and the traditional K\"ahler geometry through a series
important conjectures and theorems. In 1997, following his
program, the second author proves some of
his conjectures:\\

\noindent {\bf Theorem B} \cite{chen991}{\it The following statements
are true:
\begin{enumerate}
\item The space of K\"ahler potentials ${\cal H}$ is convex by $C^{1,1}$ geodesics.
More specifically, if $\varphi_0,\varphi_1 \in \cal H$ and
$\varphi(t) \;(0\leq t \leq 1)$ is a geodesic connecting these
two points in $\cal H,$ then the second order mixed covariant
derivatives of $\varphi(t)$ are uniformly bounded from above.
\item  ${\cal H}$ is a metric space\footnote{This is a conjecture of Donaldson\cite{Dona96}.}. In other words, the infimum of the
lengths of all possible curves between any two different points in
$\cal H$ is strictly positive.
\item If $C_1(V) < 0,$ then the extremal K\"ahler metric is unique in each
K\"ahler class.
\end{enumerate}
}
\subsection{An approximate geodesic Lemma}
In a local coordinate of $V,$ let $\omega_0 = \displaystyle \sum_{\alpha, \beta=1}^n\; {g_0}_{\alpha \overline{\beta}}
 \; d\;z_{\alpha}\;\overline{d\;z_{\beta}}$ and
 \[ \omega_{\varphi} =  \displaystyle \sum_{\alpha, \beta=1}^n\; {g}_{\alpha \overline{\beta}} \;
 d\;z_{\alpha}\;\overline{d\;z_{\beta}}, \qquad {\rm where}\qquad
 {g}_{\alpha \overline{\beta}} = {g_0}_{\alpha \overline{\beta}}
+  {{\partial^2 \varphi} \over {\partial z_{\alpha} \partial
\overline{z_{\beta}}}}.\] In this subsection,   $z_1, z_2 \cdots,
z_n$ are local coordinates in $V;$  and we always use the
following  notations: \[
 i,j,k =1,2,\cdots n,n+1 \qquad {\rm and}\qquad
\alpha,\beta, \gamma=1,2,\cdots n.\] For any path $\varphi(\cdot,
t): [0,1] \rightarrow \cal H, $ we can view  it as  a function
defined in the product manifold $V \times [0,1].\;$ Following an
idea of S. Semmes, we introduce a dummy variable $\theta$ such
that $V \times ([0,1]\times S^1)$ is a $(n+1)$ dimensional
K\"ahler manifold and $t = re(z_{n+1}).\;$ Here $S^1$ is the unit
circle. Consider the projection
\[ \begin{array} {rcl} \pi: V \times ([0,1] \times S^1) & \rightarrow  &
V\\                    (z,t,\theta) & \rightarrow & z.
\end{array}
\]
  Consider the pull back metric
$\pi^*g_0.\;$  Note that $\pi^* \omega_0$ is a degenerate
K\"ahler form  of co-rank 1 in $ V \times ([0,1]\times S^1).\;$
\begin{defi} A path $\varphi(t) (0< t < 1)$ in $\cal H$
 is a convex path if
\[
 \det\; \left( {\pi^*g_0}_{i \overline{j}}\; + \; {{\partial^2 \varphi} \over
{\partial z_i \partial \overline{z_j}}} \right)_{(n+1) (n+1)} > 0,
\qquad \;{\rm in}\;V\times (I\times S^1). \]
\end{defi}

\begin{defi}
A convex path $\varphi(t)$ in the space of K\"{a}hler metrics is
called an $\epsilon$-approximate geodesic if the following holds:
\begin{eqnarray}
 \det\; \left( {\pi^*g_0}_{i \overline{j}}\; + \; {{\partial^2 \varphi} \over
{\partial z_i \partial \overline{z_j}}} \right)_{(n+1) (n+1)}
 & = &
 ({{\p^2  \varphi} \over {\p t^2}} - {1\over 2} |\nabla {{\p  \varphi} \over {\p t}}|_{g(t)}^2)\; \det\; g(t) \nonumber\\
 &  = & \epsilon \cdot \det\; g_0
\label{eq:epsilongeodesic}
\end{eqnarray}
where $g(t)_{\alpha\overline{\beta}} =
{g_0}_{\alpha\overline{\beta}} + {{\partial^2
\varphi}\over{\partial z_{\alpha}\partial
\overline{z}_{\beta}}}\; (1\leq \alpha,\beta\leq n).\;$
\end{defi}

Also from \cite{chen991}, we have the following:
\begin{lem} \cite{chen991}(Geodesic approximation lemma):
 Suppose $\phi_1(\cdot,s), \phi_2(\cdot,s):
[0,1] \rightarrow {\cal {H}} $ are two smooth curves in ${\cal
{H}}.\;$ For $\epsilon_0$ small enough, there exist two
parameters smooth families of curves $
\varphi(\cdot,t,s,\epsilon): [0,1]\times [0,1]\times
(0,\epsilon_0]  ( 0\leq t, s \leq 1,$ and $ 0 < \epsilon \leq
\epsilon_0)\rightarrow \cal H$ (from $\phi_1(\cdot,s)$ to
$\phi_2(\cdot,s)$) such that the following properties hold:
\begin{enumerate}
\item For any fixed $s $ and $\epsilon, \varphi(\cdot, t, s,\epsilon)$ is an $\epsilon$-approximate geodesic
connecting $\varphi_1(\cdot,s)$ and $\varphi_2(\cdot,s).\;$ More
precisely, $\varphi(\cdot,t,s,\epsilon) $ solves the corresponding
Monge-Ampere equation
\begin{equation}
\det\;(\pi^* {g_0}_{i \bar j} + {{\partial^2
\varphi}\over{\partial z_{i}
\partial \overline{z}_{j}}}) = \epsilon \cdot  \det\;(g_0),
    \; {\rm in}\; V\times {\bf R}; \end{equation}
    and
\begin{equation}
\varphi(\cdot,0,s,\epsilon) =
\phi_1(\cdot,s),\;\varphi(\cdot,1,s,\epsilon) = \phi_2(\cdot,s).
\nonumber
\end{equation}
Here  $\varphi$ is independent of $Im(z_{n+1}).\; $
\item There exists a uniform constant $C$ which depends only on $\phi_1(\cdot,s), \phi_2(\cdot,s)$ such that
\[
   |\varphi| + |{{\partial \varphi}\over {\partial s}} | + |{{\partial \varphi}\over {\partial t}} | < C; \qquad  0< {{\partial^2 \varphi}\over {\partial t^2}}  < C,
\qquad {{\partial^2 \varphi}\over {\partial s^2}} < C.\]
\item For fixed $s,$ let $\epsilon \rightarrow 0,$ the $\epsilon-$ approximating
geodesic $ \varphi(\cdot, t, s,\epsilon)$ converges to the unique
geodesic between $\phi_1(\cdot,s)$ and $\phi_2(\cdot,s)$ in weak
$C^{1,1}$ topology.
\item Define energy element  along $\varphi(\cdot, t, s,\epsilon)$ by
\[E(t,s,\epsilon ) = \displaystyle \int_{V} |{{\partial \varphi}\over {\partial t}}|^2 d\;g(t,s,\epsilon), \]
where $g(t,s,\epsilon)$ is the corresponding K\"{a}hler metric
defined by the K\"ahler potentials $\varphi(t,s,\epsilon).\;$ Then
there exists a uniform constant $C$ such that
 \[
\max_{t,s}  |{{\partial \, E}\over {\partial\,t}}| \leq \epsilon
\cdot C. \] In other words, both the energy and length element
converge to a constant along each convex curve if $\epsilon
\rightarrow 0.\;$
\end{enumerate}
\end{lem}
This is a crucial lemma needed in the proof in the next subsection.

\subsection{Length of Jacobi vector field grows super-linearly}
In this subsection,  we use the same notation as in Lemma 2.3.  We
want to prove that the Jacobi vector field along any geodesic
grows sup-linearly.
 \begin{lem} Let $\varphi(\cdot, t,s,\epsilon) $ be the two parameter families of
 approximating geodesics defined as in Lemma 2.3. Let $Y(\cdot,t,s,\epsilon) = {{\p
 \varphi} \over {\p s}}$ be the deformation vector fields and $X(\cdot,t,s,\epsilon) =
 {{\p \varphi}\over {\p s}}$ the tangential vector fields along
 the  approximating geodesic. Then the second derivatives of $Y$
 along the approximating geodesic are positive:
 \[
   \nabla_X \nabla_X Y \geq 0.\]
   Note that $Y$ converges to a Jacobi vector field as $\epsilon
   \rightarrow 0.\;$ Moreover, we have
\[
    \langle Y, \nabla_X Y \rangle \geq \langle Y,Y\rangle.
   \]

 \end{lem}
 \begin{proof}
  The equation for a family of $\epsilon-$ approximate geodesics is:
  \[   ({{\p^2  \varphi} \over {\p t^2}} - {1\over 2} |\nabla {{\p  \varphi} \over {\p t}}|_{g(t)}^2)\; \det\; g(t)
   = \epsilon \cdot \det\; g_0, \qquad 0\leq \;s,\;t \;\leq 1.
\]
  Denote
\[
  X = {{\partial}\over {\partial t}},\qquad Y = {{\partial}\over {\partial
  s}},\qquad
  Y' = \nabla_X Y; \]
 and
 \[H = {{ \det\; g_0}\over {\det\; g}},\qquad
f = {{\p^2  \varphi} \over {\p t^2}} - {1\over 2} |\nabla {{\p
\varphi} \over {\p t}}|_{g(t)}^2 = \nabla_X X.\]
 Then the approximating geodesic
equation becomes
\[
f = \nabla_X X = \epsilon \cdot H.
\]
Note  that  any two-parameter family of smooth functions $F(s,t)$
can be viewed as a two parameter family of  tangent vectors at
$T_{\varphi(\cdot,t,s,\epsilon)} {\cal H}.\;$ Then, the Riemannian
metric in $\cal H$ gives the following covariant derivatives:
\[
 \nabla_X F = \nabla_{{{\partial }\over {\partial t}}} F(s,t) = {{\partial F}\over {\partial t}} - {1\over 2} \nabla_g {{\partial \varphi }\over {\partial t}}
\cdot \nabla_g {{\partial F(s,t) }\over {\partial t}}\]
and
\[
\nabla_Y F = \nabla_{{{\partial }\over {\partial s}}} F(s,t) = {{\partial F}\over {\partial s}} - {1\over 2} \nabla_g {{\partial \varphi }\over {\partial s}}
\cdot \nabla_g {{\partial F(s,t) }\over {\partial s}}.
\]
Clearly, as $\epsilon \rightarrow 0,$ $Y$ is the Jacobi vector field
along the geodesic. By definition, the length of $Y$ at $t$ is:
\[
|Y|^2(t,s) = \displaystyle \int_V \;|{{\partial \varphi}\over
{\partial s}} |^2 \; \det\; g.
\]
Then
\[
 {1\over 2} {{\partial \,}\over {\partial \,t}} |Y|^2  =
 \langle \nabla_X Y, Y\rangle = \langle \nabla_Y X, Y\rangle.
\]
Let $K(X,Y)$ denote the sectional curvature of the space of
K\"ahler metrics at point $\varphi(\cdot,t,s,\epsilon).\;$  By a
formal calulation (cf.  \cite{Ma87}, \cite{Semmes92} and
\cite{Dona96}), we have\footnote{ For any two functions $f_1, f_2$
and a K\"ahler form $\omega_{\varphi}$,  the term
$\{f_1,f_2\}_{\varphi}$ is defined to be the Possion brake of
$f_1$ and $f_2$ with respect to the sympletic form
$\omega_{\varphi}.$ }
\[ K(X,Y) = - | \{ X,Y\}_{\varphi}|_g^2 \leq 0.\]
 Therefore, we have
\[
\begin{array} {ccl}
{1\over 2} {{\partial^2 \,}\over {\partial \,t^2}} |Y|^2
 & = &  \langle \nabla_Y X, \nabla_X Y\rangle +  \langle \nabla_X\;\nabla_Y X, Y\rangle \\
 & = & |Y'|^2 - K(X,Y) + \langle \nabla_Y\;\nabla_X X, Y\rangle \\
 & \geq  & |Y'|^2 + \displaystyle \int_V \epsilon {{\partial \varphi}\over {\partial s}}
\nabla_{{{\partial} \over {\partial s}}} H\; \det g \\
& = & |Y'|^2 + \displaystyle \int_V \epsilon {{\partial \varphi}\over {\partial s}}\cdot  ( {{\partial H }\over {\partial s}}  - {1\over 2} \nabla {{\partial \varphi}\over {\partial s}} \cdot \nabla H )\; det g\\
& = & |Y'|^2 + {{\epsilon}\over {2}} \displaystyle \int_V |\nabla
{{\partial \varphi}\over {\partial s}} |^2 H  \cdot \det\;
g\;\geq |Y'|^2.
\end{array}
\]
The last equality holds since
\[{{\partial H }\over {\partial s}} = {{\partial }\over {\partial s}} ({{\det\; g_0}\over {\det\;g}}) = - \triangle_{g} {{\partial \varphi }\over {\partial s}} \cdot H
\]
and
\[
\begin{array}{lcl} & &
- \displaystyle \int_V {{\partial \varphi }\over {\partial s}}
\triangle_{g} {{\partial \varphi }\over {\partial s}} \cdot H
\;\det \;g\\
 &  = &  {1\over 2} \displaystyle \int_V |\nabla
{{\partial \varphi}\over {\partial s}} |^2 H  \cdot \det\; g +
{1\over 2} \displaystyle \int_V {{\partial \varphi }\over
{\partial s}} \nabla {{\partial \varphi }\over {\partial s}}
\nabla H \;\det \;g.
\end{array}
\]
It follows  that
\[
  {{\partial^2 \,}\over {\partial \,t^2}} |Y| \geq 0.
\]
In other words,  $ |Y(t)| (0 \leq t \leq 1)$ is a  convex function
of $ t.\; $ Since $Y(0)= 0,$  we have
\[
 {{\partial  } \over { \partial t}} |Y (t)|_{t=1}\geq {{|Y(1)| }\over { 1}}.
\]
or at time $t=1$
\begin{equation}
 \langle Y, Y' \rangle \geq \langle Y, Y\rangle.
\end{equation}
\end{proof}
 \subsection{Proof of Theorem 1.1}
In this subsection, we want to show that ${\cal H}$ is a
non-positive curved space. We follow again the notations in Lemma
2.3 and the preceding subsection.
\begin{proof} Consider a special case of Lemma 2.3 when
$\phi_1(\cdot, s) = \phi_1 $ is one point on $\cal H\;$ (instead
of a curve).  We denote this point as   $P.\; $  Let $Q =
\phi_2(\cdot,0)\in {\cal H}$ and $ R = \phi_2(\cdot,1)\in \cal
H.\;$  Furthermore, we assume that $\phi_2(\cdot,s)$ (denoted as
by $QR$) is an $\epsilon-$approximate geodesic connecting $Q$ and
$R.\;$ In other words, it satisfies the following equation:
\[
\nabla_Y Y \cdot \det\; g =
 ({{\p^2  \varphi} \over {\p t^2}}_{ss} - {1\over 2}  |\nabla {{\p  \varphi} \over {\p t}}_s|^2_{g}) \;\det\; g = \epsilon \cdot \det\; g_0.
\]
Denote $Q(s)$ the point $\phi_2(\cdot,s); $ and denote $E(s)$  the
energy of $\epsilon-$approximate geodesic from $P$ to $Q(s).\;$
As $\epsilon \rightarrow 0, E(s) \rightarrow $ a constant which,
by our normalization, is  the square of the geodesic distance from
$P$ to $Q(s).\;$ Thus it is enough to work with $E(s).\;$ Next
\[
 E(s) = \displaystyle \int_0^1 \langle X,X\rangle d\,t = \displaystyle \int_0^1 \displaystyle \int_V \;({{\partial \varphi}\over {\partial t}})^2\;
  \det\;g\; d\;t
\]
and
\[
  E(QR) = \displaystyle \int_0^1 \langle Y,Y\rangle d\,s = \displaystyle \int_0^1 \displaystyle \int_V \;({{\partial \varphi}\over {\partial s}})^2\; \det\;g\;d\;s.
\]
Thus
\[
\begin{array} {ccl}
{1\over 2} {{d \, E(s) }\over {d \,s}}  & = &
 \displaystyle \int_0^1 \langle \nabla_Y X, X\rangle d\,t   =  \displaystyle \int_0^1 (X \langle X, Y\rangle - \langle \nabla_X X, Y\rangle) d\,t \\
& = & \langle X,Y\rangle_{t=1} - \displaystyle \int_0^1 \displaystyle \int_V {{\partial \varphi}\over {\partial s}} \cdot \epsilon\; H \;\det g \; d t\\
& = &  \langle X,Y\rangle_{t=1} - \epsilon \cdot \displaystyle
\int_0^1 \displaystyle \int_V {{\partial \varphi}\over {\partial
s}}\; \det\;g_0\; d\,t.
\end{array}
\]
Now the second derivatives:
\[
\begin{array} {ccl}
{1\over 2} {{d^2 \, E(s) }\over {d \,s^2}}  & = &  {d\over {d\,s}}  \langle X,Y\rangle_{t=1} - \epsilon \cdot \displaystyle \int_0^1 \displaystyle \int_V {{\partial^2 \varphi}\over {\partial s^2}}\; det\;g_0\; d\,t\\
&\geq & \langle Y',Y\rangle_{t=1} + \langle X, \nabla_Y Y\rangle_{t=1} - C\;\epsilon \displaystyle \int_V \det\;g_0\\
& \geq  & \langle Y,Y\rangle_{t=1}+ \displaystyle \int_V {{\partial \varphi}\over {\partial t}} \;\epsilon \cdot H \cdot \det g - C\;\epsilon \;\displaystyle \int_V \det\;g_0 \\

& \geq & \langle Y,Y\rangle_{t=1} + \displaystyle \int_V {{\partial \varphi}\over {\partial t}} \;\epsilon \cdot \det\;g_0 - C\;\epsilon \;\displaystyle \int_V \det\;g_0\\
& \geq & E(QR)- C\cdot \epsilon \cdot vol(V).
\end{array}
\]
Here we have used the inequality (2.6) in the second inequality
from the top. And $E(QR)$ denotes the energy of the path
$\phi_2(\cdot,s).\;$ For the energy elements of curves, the
following inequality holds
\[
  E(s) \leq (1-s) E(0) + s E(1) - s(1-s) (E(QR) - C\cdot\epsilon \cdot vol(V))
\]
Now fix $s$\footnote{Actually, using successive subdivision one
sees that knowing the inequality (2.7) holds for $s= {1\over 2}$
is enough to prove it  for all $0 < \lambda < 1,$ cf.
\cite{SchoenKore93}.}, let $\epsilon \rightarrow 0,$ each energy
element of a path approaches the square of the length of that
path. Thus the above inequality reduces to:
\begin{equation}
  |PQ(s)|^2 \leq (1-s) |PQ|^2 + s |PR|^2 - s(1-s) |QR|^2.
\end{equation}
Thus the space of K\"ahler  metrics satisfies the defining
inequality for non-positive curved space and hence it is a
non-positive space. Here $|PQ(s)|$ represents the distance from
$P$ to $Q(s)$; $|PQ|$ represents the distance from $P$ to $Q$;
$|PR|$ represents the distance from $P$ to $R$; and $|QR|$
represents the distance from $Q$ to $R.\;$
\end{proof}

Next we  prove Theorem 1.3.
\begin{proof} Let $\varphi_0$ and $\varphi_1$ be two points in
$\cal H$ with distance $l > 0.\;$ Suppose that $\varphi(\cdot, t)
(0 \leq t \leq 1) $ is a $C^{1,1}$ geodesic which connects these
two points in $\cal H.\; $  Let  $\varphi_i(t)(0 \leq t \leq 1) $
be an arbitrary  family ($i=1,2,\cdots n\cdots$) of curves between
$\varphi_0$ and $\varphi_1 $ with length $l_i \geq l >0.\;$ Next
we assume that this is a distance minimizing sequence of curves.
In other words, \[ \displaystyle \lim_{i \rightarrow \infty} l_i =
l.\]
 Then,
 we need to show that $\varphi_i(\cdot, t) (0 \leq t \leq 1) $ converges to
$\varphi(\cdot, t) (0 \leq t \leq 1) $ in some reasonable
topology. For convenience, we assume that  every curve involved
has been parameterized proportional to the arc-length. Then we
only need to show that for each fixed $s> 0$, we have
\[
\displaystyle \lim_{i \rightarrow \infty}\; d(\varphi_i(\cdot,s),
\varphi(\cdot, s)) \rightarrow 0. \]
 Since $\cal H$ is a non-positive curved space, we have (comparing with the Euclidean
space):
\[
  d(\varphi_i(\cdot, s), \varphi(\cdot, s) \leq \sqrt{{{l_i^2 - l^2} \over {4}}}\rightarrow 0.
\]
Theorem 1.3 is then proved.
\end{proof}
\section{The  gradient flow of the K energy}
In this section, we  prove Theorem 1.4.

\begin{proof}
Let $\varphi(\cdot, t): [0,1] \rightarrow \cal H$ be any smooth
curve in $\cal H.\;$ Suppose that $\varphi(\cdot,t,s)$ is the
image of this curve under the gradient flow after time $s.\;$
Recall \[ {{\p \varphi}\over {\p s}} = R(\varphi) - \underline{R}.
\]
 Denote  $g(s,t)$ as the K\"ahler metric associated with the
K\"ahler potentials $\varphi(s,t).\;$ Use $\triangle $ to denote
the complex Laplacian operator of metric $g(s,t).\;$  Following a
calculation in \cite{calabi82}, we have
\[
   {{\partial R} \over {\partial t}} = - D^* D \;{{\partial \varphi} \over {\partial t}} + \displaystyle \sum_{\alpha=1}^n\;({{\partial \varphi} \over {\partial t}})\varphi^{\alpha} R_{,\alpha}
\]
and \[
 {{\partial } \over {\partial t}}\det\; g = \triangle \;
{{\partial \varphi} \over {\partial t}} \; \det\; g.\]
 Recall that the
energy of the path $\varphi(\cdot, t, s)$ (at time $s$ fixed) is:
\[
  E(s) = \int_0^1 \int_V\; \left({{\partial \varphi} \over {\partial t}}\right)^2\; \det g\; d\, t.
\]
Under the gradient flow  (\ref{eq:calabiflow}), we have
\[
\begin{array} {lcl}
{{d E}\over {d\,s}} & = &  \int_0^1 \int_V  2 {{\partial \varphi}
\over {\partial t}} {{\partial^2 \varphi} \over {\partial t
\partial s}} \;\det \; g \;d\,t
  + \int_0^1 \int_V \;  \left({{\partial \varphi} \over {\partial t}}\right)^2 \;\triangle {{\partial \varphi} \over {\partial s}} \det\; g \; d\,t\\
  &= & \int_0^1 \int_V \; 2 {{\partial \varphi} \over {\partial t}} {{\partial R} \over {\partial t}} \;\det \; g \;d\,t
  - \int_0^1 \int_V \; 2 {{\partial \varphi} \over {\partial t}}  ({{\partial \varphi} \over {\partial t}})^{\alpha}   ({{\partial \varphi} \over {\partial s}})_{\alpha}  det\; g \; d\,t \\
& = & \int_0^1 \int_V \; 2 {{\partial \varphi} \over {\partial
t}} (  - D^* D \;{{\partial \varphi} \over {\partial t}} +
\varphi^{,\alpha} R_{,\alpha}) \;\det \; g \;d\,t \\
& & \qquad \qquad  -  \int_0^1 \int_V \;2 {{\partial \varphi}
\over {\partial t}} ( {{\partial \varphi} \over {\partial
t}})^{\alpha}
    R _{\alpha}  \det\; g \; d\,t \\
& = & -\int_0^1 \;\int_V | D {{\partial \varphi} \over {\partial
t}}|^2_g \; \det\; g \;d\;t.
\end{array}
\]
It follows that
\[
{{d L}\over {d\,s}} =  - \;\int_0^1 \left(\int_V |D {{\partial
\varphi} \over {\partial t}}|_{\varphi(s,t)}^2 \; d \,g(s,t)
\cdot \sqrt{ \int_V |{{\partial \varphi} \over {\partial t}}|^2\;
d \,g(s,t)}^{-{1\over 2}} \right)\; d\,t,
\]
where $L(s)$ is the length of the evolved curve at time $s >
0.\;$ From this formula, if the length of  a smooth curve is not
decreasing, then
\[
\int_0^1 \left(\int_V |D {{\partial \varphi} \over {\partial
t}}|_{\varphi(s,t)}^2 \; \det\;g(s,t) \sqrt{ \int_V |{{\partial
\varphi} \over {\partial t}}|^2\; d \,g(s,t)}^{-{1\over
2}}\right)\;d\;t =  0.\] It follows that
\[
   \int_0^1 \int_V |D {{\partial \varphi} \over {\partial t}}|_{\varphi(s,t)}^2
= 0
\]
or
\[
  \left({{\partial \varphi} \over {\partial t}}\right)_{,\alpha \beta}  \equiv 0, \qquad \forall \;\alpha,\beta = 1,2,\cdots n;\; \forall \; t\in [0,1].
\]
In other words, the curve $\varphi(t) (0 \leq t \leq 1) $ is
either trivial (depending only on $t$) or it represents a family
of holomorphic transformation.  Theorem 1.4 is then proved.
\end{proof}

Next we give a proof of the first part of Theorem 1.5.
\begin{proof}
 For any $\varphi_0,\varphi_1 \in \cal H, $ consider
the space of all smooth curves which connect $\varphi_0$ with
$\varphi_1.\;$ We denote it by $ {\cal
L}(\varphi_0,\varphi_1).\;$ For any curve $c \in {\cal
L}(\varphi_0,\varphi_1),$ we denote its length by $L(c).\;$ Then
the distance between the two points $\varphi_0$ and $\varphi_1$
can be defined as \[ d(\varphi_0,\varphi_1) = \displaystyle
\inf_{c \in {\cal L}(\varphi_0,\varphi_1)} L(c).
\]
We also define a map in $\cal H$ via the gradient flow (1.1): for
a fixed time $s,$ and for any $\varphi \in \cal H,$ we define
that the image of $\varphi$ under the map $\pi_s$ is the image of
$\varphi$ along the gradient flow after time $s>0,\;$ provided
the gradient flow initiated at $\varphi$  does exist for time
$s>0.\;$ It is clear that for any $\varphi$, the map is defined
for small $s> 0.\;$ However, for a fixed $s>0$, $\pi_s$ is not
necessarily defined for all $\varphi \in H $ since we don't know
the global existence of the gradient flow.

On the other hand, if the gradient flow exists for all the time
for any smooth initial metric, then this induces a well defined
map from ${\cal L}(\varphi_0,\varphi_1)$ to ${\cal L}(\pi_s
(\varphi_0),\pi_s(\varphi_1))$ for any $s>0.\;$ Since the length
of any smooth curve in $\cal H$ is decreased under the gradient
flow, we have

\[
\displaystyle \inf_{c \in {\cal
L}(\pi_s(\varphi_0),\pi_s(\varphi_1))} L(c) \leq \displaystyle
\inf_{c \in {\cal L}(\varphi_0,\varphi_1)} L(c), \qquad \forall\;
s > 0 .
\]
Thus,
\[ d(\pi_s(\varphi_0), \pi_s(\varphi_1)) \leq
d(\varphi_0,\varphi_1), \qquad \forall\; s > 0.
\]
\end{proof}
 Before proving the second part of Theorem 1.5, we need to use a
theorem in \cite{chen991} where an explicit formula for the first
 derivatives of the distance function in $\cal H$  is given.
For the convenience of the readers, we will include this theorem
here. The second part of Theorem 1.5 is essentially a corollary of
this theorem.

\begin{theo}\cite{chen991}  For any two K\"{a}hler potentials $\varphi_0, \varphi_1,$
the distance function $d(\varphi_0,\varphi_1)$ is at least $C^1.\;$ More
specifically, if $\varphi_0, \varphi_1$ move along two curves $\varphi_0(s),
$ and $\varphi_1(s)$ respectively, and if we denote the distance between
$\varphi_1(s)$ and $\varphi_2(s)$ is $L(s),$ then
\[ \begin{array}{lcl} {{d\; L(s)} \over {d\,s}}\mid_{s=0} &  =  & \langle X,Y_1\rangle |X|^{-{1\over 2}}
\mid_{t=1}- \langle X,Y_0\rangle |X|^{-{1\over 2}} \mid_{t=0} \\
& = &  \displaystyle \int_V {{\partial \varphi_1}\over {\partial
s}} \;{{\partial \, \varphi}\over {\partial\, t}}\; d\, g(s)
\cdot \{\displaystyle \int_V |{{\partial \varphi}\over {\partial
t}}|^2 d\,g(s)\}^{-{1\over 2}}\mid_{t=1} \\
& & \qquad -  \displaystyle \int_V {{\partial \varphi_0}\over
{\partial s}} \;{{\partial \, \varphi}\over {\partial\, t}}\; d\,
g(s) \cdot \{\displaystyle \int_V |{{\partial \varphi}\over
{\partial t}}|^2 d\,g(s)\}^{-{1\over 2}}\mid_{t=0}.
\end{array}
\]
Here $ \varphi(t) (0 \leq t \leq 1)$ denotes the $C^{1,1}$
geodesic connecting the two metrics $\varphi_0$ and
$\varphi_1;\;$ and $X = {{\partial \varphi}\over {\partial t}}
\in T_{\varphi(t)} \cal H$ and $Y_i = {{\partial \varphi_i}\over
{\partial s}} \in T_{\varphi_i} {\cal H}\;(i=0,1).$
\end{theo}
Now we complete the proof of Theorem 1.5.
\begin{proof} If the gradient flow (\ref{eq:calabiflow}) exists for
all the time, then it is straightforward to show that flow
(\ref{eq:calabiflow}) decreases the distance between any two
points in $\cal H $ unless they are connected by a holomorphic
transformation. Thus, we only deal with the case when the K
energy is weakly convex. By definition, for any curve $\varphi(t)
\in \cal H,$ the K energy is defined as
\[
{{d\, M (\varphi(t))}\over {d\,t}} = - \displaystyle \int_V \;
{{\partial \varphi}\over {\partial t}}\; (R - \underline{R}) \;
\det g.
\]
Along a $C^{1,1}$ geodesic, the second derivative of the K energy
is convex in the weak sense that
\[
{{d^2\, M (\varphi(t))}\over {d\,t^2}} \geq 0.
\]
In particular, we have
\begin{equation}
  {{d\,M(\varphi(t))}\over {d\,t}} \mid_{t=1} \geq  {{d\,M(\varphi(t))}\over {d\,t}}
  \mid_{t=0}.
\label{eq:mabuchiweak}
\end{equation}

 Suppose that $\varphi(t) (0\leq t \leq 1) $ is the unique $C^{1,1}$
geodesic which connects $\varphi_1$ and $\varphi_2,$ and suppose
it is parameterized proportional to arc length. If we flow
$\varphi_1$ and $\varphi_2$ by the gradient flow
(\ref{eq:calabiflow}), we have

\[
  {{\partial \varphi_1} \over {\partial s} } = R(\varphi_1(s)) - \underline{R} \qquad {\rm and}\;  {{\partial \varphi_2} \over {\partial s} } = R(\varphi_2(s)) - \underline{R}.
\]
Plugging this into the corresponding formula in Theorem 3.1, we
have
\[
\begin{array}{lcl}  {{d\; L(s)} \over {d\,s}} & = & \{\displaystyle \int_V
|{{\partial \varphi}\over {\partial t}}|^2 d\,g(s)\}^{-{1\over
2}} \\
& & \qquad \cdot  \left( \displaystyle \int_V (R(\varphi_2) -
\underline{R}) \;{{\partial \, \varphi}\over {\partial\, t}}\;
d\, g(s)\mid_{t=1} -   \displaystyle \int_V
(R(\varphi_1) - \underline{R}) \;{{\partial \, \varphi}\over {\partial\, t}}\; d\, g(s)\mid_{t=0}\right) \\
& = & -\{\displaystyle \int_V
|{{\partial \varphi}\over {\partial t}}|^2 d\,g(s)\}^{-{1\over 2}} \cdot \left( {{d\;M}\over {d \,t}} \mid_{t=1} - {{d\;M}\over {d \,t}} \mid_{t=0}\right) \\
& = & - \{\displaystyle \int_V
|{{\partial \varphi}\over {\partial t}}|^2 d\,g(s)\}^{-{1\over 2}} \cdot \int_{t=0}^1 {{d^2\;M}\over {d \,t^2}} \; d\;t \leq 0.
\end{array}
\]
\end{proof}


\section{Some further corollaries, remarks and the relationship with stability}

\begin{cor} If all geodesics are smooth, then the extrmal K\"ahler
metric is unique up to some holomorphic automorphisms.
\end{cor}
\begin{proof}
  Suppose that there exist two extremal
K\"ahler metrics in a fixed K\"ahler class. It was proved in
\cite{calabi85} that any extremal K\"ahler metric must be
symmetric with respect to a maximal compact sub-group. Without
loss of generality, one may assume that both metrics are
symmetric under the same maximal compact sub-group. Then every
metric in the geodesic which connects this two extremal metrics
must also have the same symmetry group (via Maximum principle).
If the scalar curvature is constant, then an argument of Donaldson
\cite{Dona96} on the convexity of the K energy implies that
 the extremal metric must be unique. If the scalar curvature is not a constant,
then the gradient vector field of the scalar curvature is a
holomorphic vector field and it is unique in each K\"ahler class
once the maximal compact sub-group is fixed. In particular, the
gradient flow  (\ref{eq:calabiflow}) restricted to the two
extremal metrics induces the exact same holomorphic
transformation. Thus the distance of these two extremal K\"ahler
metrics under the gradient flow is unchanged. Suppose $\varphi_1,
\varphi_2 $ are the two extremal metrics and $\varphi(\cdot,s)$
is the unique geodesic connecting them.  Since the distance of
$\varphi_1$ and $\varphi_2$ is not decreased under the gradient
flow,  by Theorem 1.3, the path $\varphi(\cdot,s)$ must either be
totally trivial or represent a holomorphic transformations.
\end{proof}
\begin{rem} Donaldson in 1997 gave a proof to this corollary in the case of constant scalar
curvature metric;  and suggested to the second author  that a
modified proof of his  works for general extremal K\"ahler
metrics. In this paper, we presented a new proof.
\end{rem}
\begin{rem}
For the uniqueness of the extremal K\"{a}hler metric, the  known
results are as follows:
  1)in 1950s, the first author showed the uniqueness of
K\"{a}hler-Einstein metric if $C_1 \leq 0.\;$
  2)in 1987, T. Mabuchi and S. Bando~\cite{Bando87} showed the uniqueness
of K\"{a}hler-Einstein metrics up to holomorphic transformation
if the first Chern class is positive. In \cite{chen991}, the 2nd
author proved that the constant scalar curvature metric is unique
in each K\"ahler class if $C_1(V) < 0.\;$ The problem for the
general case is still open. However, the second author
\cite{chen943}
 had examples of non-uniqueness of some degenerated
extremal K\"ahler metrics in $S^2.\;$
\end{rem}

S. T. Yau predicted in \cite{Yau92} that the existence of
K\"ahler-Einstein metrics is related to the stability in the
sense of Hilbert Schemes and Geometric invariant theory. His
conjecture should be extended to include the case of extremal
K\"ahler metrics.  From Theorem 1.3, we observe some kind of link,
perhaps still a bit mysterious, between the the existence of
extremal metrics and ``stability" of the infinite dimensional
space $\cal H$  in some sense. At least formally, it fits nicely
to the general picture  Yau's conjecture describes. The following
paragraph is essentially speculative in the effort to explaining
this point. If we are willing to put aside the regularity issue,
then Theorem 1.3 implies that the gradient flow of the K energy
is a distance contracting flow in $\cal H$. In this infinite
dimensional path length space $\cal H$, we choose a large enough
ball, which hopefully contains any possible candidates for
extremal K\"ahler metrics. Now flow the entire ball by this
gradient flow, if global solution of the gradient flow always
exist for all smooth initial metric, then the contracting nature
of the flow will shrink the size of the ball. At the end, the ball
shall be contracted to a point,  and the limit point must be an
extremal K\"ahler metric we are looking for. However, this formal
picture is not quite complete. A dichotomy can possibly taken
place: As the size of the ball shrinks, the ball may also be
drifted away to infinity. In the first possibility when the ball
stays in a finite domain,  the infinite dimensional manifold is
considered "stable"  in some sense and we arrive at the unique
extremal K\"ahler metric in the limit of the flow. In the second
case when the ball drifts to the infinity, then the infinite
dimensional space is considered as "unstable" in some sense, and
the gradient flow converges to an extremal K\"ahler metric in a
different K\"ahler manifold.

 \end{document}